\pdfoutput=1
\documentclass[12pt]{amsart}
\usepackage{amssymb,mathtools,pgfplots}
\usepackage{graphicx,float,array,booktabs,tabularx,longtable}
\definecolor{hot}{RGB}{65,105,225}
\definecolor{coderule}{RGB}{190,190,184}
\definecolor{codecomment}{RGB}{90,100,90}
\definecolor{codestring}{RGB}{120,75,45}
\definecolor{sagegreen}{RGB}{0,100,0}

\usepackage[colorlinks, allcolors=hot]{hyperref}
\usepackage[backend=biber,style=numeric,sorting=nyt,maxbibnames=99]{biblatex}
\usepackage[margin=1in]{geometry}
\usepackage{algorithm}
\usepackage[noend]{algpseudocode}
\newcommand{\figcaptiongap}{\par\vspace{0.6\baselineskip}}

\pgfplotsset{compat=1.18}

\usepackage{listings}

\algrenewcommand\algorithmicrequire{\textbf{Input:}}
\algrenewcommand\algorithmicensure{\textbf{Output:}}
\algrenewcommand\algorithmiccomment[1]{\hfill{\(\triangleright\) #1}}

\newtheorem{theorem}{Theorem}

\theoremstyle{definition}
\newtheorem{definition}[theorem]{Definition}

\newcommand{\T}{\mathbb{T}}

\title[KnotMosaics for SageMath]{The \texttt{KnotMosaics} Package for SageMath}

\author{Mary Y. Deng}
\address[M. Y. Deng]{Dept. of Mathematics\\
University of Washington\\
Seattle, WA 98195\\
USA}
\email{marydeng@uw.edu}

\author{Allison K. Henrich}
\address[A. K. Henrich]{Mathematics Department\\
Seattle University\\
Seattle, WA 98122\\
USA}
\email{henricha@seattleu.edu}
\urladdr{https://allisonhenrich.com}

\author{Sean H. Kawano}
\address[S. H. Kawano]{Dept. of Mathematics\\
University of Washington\\
Seattle, WA 98195\\
USA}
\email{shkawano@uw.edu}

\author{Andrew R. Tawfeek}
\address[A.~R.~Tawfeek]{Dept. of Mathematics\\
         University of Washington\\
         Seattle, WA 98195\\
         USA}
\email{atawfeek@uw.edu}
\urladdr{https://atawfeek.com/}


\begin{document}

\begin{abstract}
We introduce \texttt{KnotMosaics}, a \emph{SageMath} package for constructing, visualizing, and analyzing knot mosaic diagrams. The package represents an $n$-mosaic as a matrix of standard tile labels and implements the local connectivity rules needed to validate mosaics, trace strands and components, compute planar diagram codes, generate random examples, and construct rational tangle mosaics. The planar diagram interface connects the mosaic representation to existing knot and link software, enabling computations such as Jones polynomials and knot Floer homology checks. We describe the package design, its main algorithms, and representative examples that illustrate how \texttt{KnotMosaics} can support computational exploration in knot mosaic theory.
\end{abstract}

\maketitle


\section{Introduction}

Knot mosaic diagrams, introduced by Lomonaco and Kauffman~\cite{lomonaco2008quantum} in the context of quantum knot systems, provide a combinatorial framework for representing knots and links as $n \times n$ grids of tiles drawn from a standard set $\T = \{T_0, \dots, T_{10}\}$. The tiles are required to be \emph{suitably connected}, meaning that each connection point of a tile meets a connection point on an adjacent tile. Kuriya and Shehab~\cite{kuriya2014lomonaco} proved (Theorem~\ref{lc-conj}) that this mosaic framework faithfully captures topological knot equivalence: two tame knots are of the same type if and only if any mosaic representatives are related by a finite sequence of mosaic Reidemeister moves.

\begin{figure}[ht] 
\centering
\includegraphics[width=.9\linewidth]{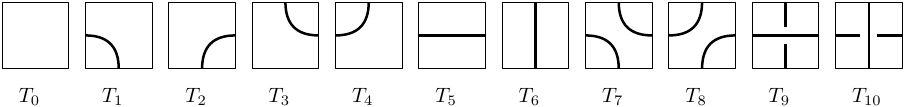} 
\figcaptiongap
\caption{The eleven standard mosaic tiles.}\label{tile-fig}
\end{figure}

Since their introduction, knot mosaics have attracted significant interest in combinatorial knot theory. Researchers have studied the \emph{mosaic number} of a knot---the minimum $n$ such that it can be represented on an $n$-mosaic---establishing bounds and tabulations for various knot families~\cite{02d_Lee2014-3-5, 02d_Ludwig2013, 02d_Ludwig2018, 02d_Dye2014}. The enumeration of knot $n$-mosaics and upper bounds for toroidal mosaics have also been investigated~\cite{02d_Lee2014-4-5, carlisle2013upper}, as have extensions to virtual knots~\cite{ganzell2020virtual}. Additionally, the authors have recently studied recursive mosaics as a means of representing wild knots~\cite{recursivemosaics}.

Despite this activity, there has been no dedicated software package for working with knot mosaics. The \texttt{KnotMosaics} package for \emph{SageMath}~\cite{sagemath} is designed to fill this gap by keeping the mathematical data model close to the literature: a mosaic is stored as a matrix of tile labels, and algorithms proceed by local checks and strand tracing. The package is intended both as a research tool for generating and testing examples and as a bridge between the mosaic literature and the broader computational knot theory ecosystem. Its main capabilities include:
\begin{enumerate}
\item constructing and visualizing knot mosaics from tile matrices;
\item verifying suitable connectivity and computing combinatorial invariants (crossing number, number of components);
\item tracing strands through mosaics and computing planar diagram codes for interfacing with knot and link software;
\item generating random mosaics subject to constraints on crossings and components;
\item constructing rational tangle mosaics and performing tangle joins; and
\item zooming mosaics by replacing each tile with an isotopy-equivalent $3 \times 3$ block, which facilitates crossing-to-crossing traversal.
\end{enumerate}

\subsection*{Software availability}
The development version of \texttt{KnotMosaics} is available at \url{https://github.com/andrew-tawfeek/knot_mosaics}. The SageMath integration described here has also been submitted for upstream review as SageMath pull request~\#42157, available at \url{https://github.com/sagemath/sage/pull/42157}. The repository accompanying this article includes a snapshot of the SageMath PR code, a pure-Python companion implementation used to generate the manuscript figures, and a test suite.

\subsection*{Outline}
In Section~\ref{sec:background}, we review the definition and basic theory of knot mosaics. In Section~\ref{sec:package}, we describe the package architecture and the algorithms used for traversal, planar diagram codes, random generation, and unknot detection. In Section~\ref{sec:examples}, we demonstrate the package through worked examples. In Section~\ref{sec:applications}, we describe computational applications and directions for future development.

\section{Background for Knot Mosaics}\label{sec:background}

Knot mosaic diagrams were first introduced by Kauffman and Lomonaco in 2008~\cite{lomonaco2008quantum} toward the development of quantum knot systems. Throughout this article, ``knot mosaic'' may refer to either a knot mosaic or a link mosaic when the distinction is not important. A mosaic diagram is an $n \times n$ grid of standard tiles satisfying a local compatibility condition called \textbf{suitable connectivity}: each connection point of a tile must meet a connection point on the adjacent tile.

\begin{definition}[Mosaic] \label{mosaicdef}
    Denote by $\mathbb{T}=\{T_0,\dots,T_{10}\}$ the standard set of 11 tiles, shown in Figure \ref{tile-fig}. For a positive integer $n$, we define the following:
    \begin{itemize}
        \item An {\bf $n$-mosaic} $M$ is an $n \times n$ matrix with entries from $\mathbb{T}$. The set of $n$-mosaics is denoted by $\mathbb{M}^{(n)}$.
        \item A {\bf knot $n$-mosaic} is an $n$-mosaic in which all tiles are suitably connected. We let $\mathbb{K}^{(n)}$ denote the subset of $\mathbb{M}^{(n)}$ of all knot $n$-mosaics.
    \end{itemize}
\end{definition}

Just as with ordinary knot diagrams, there is a set of Reidemeister moves that relate mosaics representing equivalent topological knots. This mosaic-based set of moves, however, is more extensive than the set of classical Reidemeister moves since planar isotopies must be exhaustively described in the mosaic setting, and natural operations that expand or contract the size of the mosaic must be included. See \cite{lomonaco2008quantum} for details. Kauffman and Lomonaco conjectured that their set of mosaic moves completely captured topological knot equivalence, and in 2014, Kuriya and Shehab proved it \cite{kuriya2014lomonaco}. To be precise, they proved the following. 

\begin{theorem}[Lomonaco-Kauffman Conjecture, \cite{kuriya2014lomonaco}] \label{lc-conj}
    Let $k_1$ and $k_2$ be two tame knots (or links), and let $K_1$ and $K_2$ be two arbitrarily chosen mosaic representatives of $k_1$ and $k_2$, respectively. Then $k_1$ and $k_2$ are of the same knot type if and only if the representative mosaics $K_1$ and $K_2$ are related by a finite sequence of mosaic Reidemeister moves.
\end{theorem}

Once mosaics were established as an alternative representation of topological knots and links, researchers began studying their combinatorial and geometric properties. One central invariant is the {\bf mosaic number} of a knot, the minimum number $n$ such that the knot can be represented as a knot $n$-mosaic. See, for instance, \cite{02d_Dye2014, 02d_Ludwig2018, 02d_Ludwig2013, 02d_Lee2014-3-5, 02d_Lee2014-4-5}. Other work has studied toroidal mosaics, virtual mosaics, and related geometric questions~\cite{carlisle2013upper, ganzell2020virtual}.

\section{Package Design and Algorithms}\label{sec:package}

The package is organized around a direct translation of the mathematical definition into code. The \texttt{MosaicTile} class stores the local data associated to a single standard tile $T_i \in \T$: its connection directions, strand pairings, crossing status, and local rendering instructions. The \texttt{Mosaic} class stores an $n \times n$ matrix of tile labels and recovers global information by applying these local rules across adjacent entries. This separation makes the implementation easy to audit: suitable connectivity, strand tracing, and conversion to other encodings all reduce to repeated calls to tile-level connectivity data.

\subsection{Data model}

An $n$-mosaic is instantiated from a matrix or list of lists:

\begin{lstlisting}
(*@\sageprompt@*) M = Mosaic([[0,2,1,0],[2,9,10,1],[3,10,9,4],[0,3,4,0]]); M
(*@\qquad \ @*) Mosaic of dimension 4
\end{lstlisting}

The below table gives a user-facing map of the package. It is not intended as a complete reference; instead, it organizes the routines used throughout the examples below by the mathematical task they support. Note, the method \texttt{number\_of\_crossings()} below counts crossing tiles in the mosaic diagram, \textit{not} the crossing number of the represented knot or link.

\begin{table}[H]
\centering
\footnotesize
\setlength{\tabcolsep}{3pt}
\renewcommand{\arraystretch}{1.05}
\begin{tabularx}{\linewidth}{@{}>{\raggedright\arraybackslash}p{0.17\linewidth}>{\raggedright\arraybackslash}p{0.37\linewidth}>{\raggedright\arraybackslash}X@{}}
\toprule
\textbf{Workflow} & \textbf{Representative routines} & \textbf{Role}\\
\midrule
Data model
& \texttt{Mosaic(data)}, \texttt{MosaicTile(i)}, \texttt{matrix()}
& Store square arrays of tile labels and local tile data.\\

Validation and counts
& \texttt{is\_suitably\_connected()}, \texttt{find\_crossings()}, \texttt{number\_of\_crossings()}, \texttt{number\_of\_components()}
& Check local matching conditions and compute diagram-level statistics.\\

Strand traversal
& \texttt{exit\_path()}, \texttt{strand\_of()}, \texttt{strands()}, \texttt{walk()}, \texttt{shift()}
& Trace components and move from one crossing to the next.\\

Visualization and moves
& \texttt{show()}, \texttt{zoom()}, \texttt{flip()}
& Draw mosaics, zoom tiles to $3 \times 3$ blocks, and apply basic transformations.\\

Interoperability
& \texttt{pd\_code()}, \texttt{link()}, \texttt{is\_unknot()}
& Export planar diagram data and access external knot/link computations.\\

Constructors
& \texttt{random\_mosaic()}, \texttt{rational\_tangle()}, \texttt{tangle\_join()}
& Generate constrained random mosaics and rational tangle examples.\\
\bottomrule
\end{tabularx}
\end{table}

\subsection{Traversal and planar diagram codes}

Traversal routines are the computational heart of the package. Given a starting tile and entry direction, \texttt{exit\_path} follows the local strand through the tile and identifies the adjacent tile and entry direction. Iterating this local rule gives \texttt{strand\_of}, which traces an entire component, and \texttt{walk}, which starts at a crossing and continues until it reaches the next crossing. The \texttt{zoom} operation replaces each tile by an isotopy-equivalent $3 \times 3$ block so that adjacent crossings are separated by non-crossing tiles before crossing-to-crossing walks are performed.

The planar diagram code construction is also traversal-based. Algorithm~\ref{alg:pd-code} keeps track of every still-unwalked local strand in the mosaic. Each component is walked once, arcs are labeled between consecutive crossing visits, and the two visits at each crossing determine the final planar diagram tuple.

\begin{algorithm}[H]
\caption{Planar diagram code from a mosaic}\label{alg:pd-code}
\begin{algorithmic}[1]
\Require a suitably connected mosaic $M$ with crossing set $C$
\Ensure a planar diagram code for $M$
\State $C \gets$ sorted list of crossing tiles in $M$
\State $R \gets$ set of all local strands in nonempty tiles
\While{$R \ne \emptyset$}
    \State $s \gets$ first tile with an unwalked strand
    \State walk one component from $s$
    \State delete the component's local strands from $R$
    \State label the arcs between consecutive crossing visits
    \State record crossing visits and over/under data
\EndWhile
\For{each crossing $c \in C$}
    \State order the four incident labels counterclockwise
    \State rotate the tuple to begin with the incoming under-strand label
\EndFor
\State \Return the list of crossing tuples
\end{algorithmic}
\end{algorithm}

Unlike the older \texttt{oriented\_gauss\_code} routine, the planar diagram code construction can record multiple components, and therefore applies to links as well as knots when each component has at least one crossing. Split crossing-free components require separate handling, as they are not represented by ordinary planar diagram codes.

The \texttt{pd\_code} method provides a bridge between the mosaic representation and existing knot and link software. It labels the arcs between crossing visits and records, at each crossing, the four incident arc labels in counterclockwise order beginning with the incoming under-strand arc. The resulting planar diagram can then be passed to SageMath's \texttt{Link} class, giving access to computations such as those in the table below.

\begingroup
\footnotesize
\setlength{\tabcolsep}{4pt}
\renewcommand{\arraystretch}{1.15}
\setlength{\LTleft}{0pt}
\setlength{\LTright}{0pt}
\begin{longtable}{@{}>{\raggedright\arraybackslash}p{0.32\linewidth}>{\raggedright\arraybackslash}p{0.24\linewidth}>{\raggedright\arraybackslash}p{0.36\linewidth}@{}}
\toprule
\textbf{Function} & \textbf{Description} & \textbf{Usage}\\
\midrule
\endfirsthead
\toprule
\textbf{Function} & \textbf{Description} & \textbf{Usage}\\
\midrule
\endhead
\midrule
\multicolumn{3}{r}{\footnotesize Continued on next page}\\
\endfoot
\bottomrule
\endlastfoot
\texttt{Link.jones\_polynomial()}
& Computes the Jones polynomial of the represented link.
& {\ttfamily Link(M.pd\_code()).\linebreak jones\_polynomial()}\\

\texttt{Link.alexander\_polynomial()}
& Computes the Alexander polynomial.
& {\ttfamily Link(M.pd\_code()).\linebreak alexander\_polynomial()}\\

\texttt{Link.determinant()}
& Computes the link determinant, an integer invariant derived from the Alexander polynomial.
& {\ttfamily Link(M.pd\_code()).\linebreak determinant()}\\

\texttt{Link.signature()}
& Computes the classical signature from a Seifert matrix.
& {\ttfamily Link(M.pd\_code()).\linebreak signature()}\\

\texttt{Link.seifert\_matrix()}
& Constructs a Seifert matrix for a Seifert surface of the link.
& {\ttfamily Link(M.pd\_code()).\linebreak seifert\_matrix()}\\

\texttt{Link.khovanov\_homology()}
& Computes Khovanov homology, a categorification of the Jones polynomial.
& {\ttfamily Link(M.pd\_code()).\linebreak khovanov\_homology()}\\

\texttt{Link.get\_knotinfo()}
& Looks up matching knot or link data in the KnotInfo database when available.
& {\ttfamily Link(M.pd\_code()).\linebreak get\_knotinfo()}\\
\end{longtable}
\endgroup

\subsection{Random generation and unknot filtering}

The \texttt{random\_mosaic} function builds suitably connected mosaics incrementally. Tiles are placed in row-major order, and at each position $(i,j)$ the \texttt{potential\_tiles} method computes the set of valid tile choices by checking the connection requirements imposed by the already-placed neighbors above and to the left, together with the boundary constraints on all four sides. A tile is then chosen uniformly at random from this valid set, as summarized in Algorithm~\ref{alg:random-mosaic}. This greedy approach avoids the exponential cost of generating and filtering all $11^{n^2}$ possible $n$-mosaics, and the output is verified by calling \texttt{is\_suitably\_connected()}.

\begin{algorithm}[H]
\caption{Constrained random mosaic generation}\label{alg:random-mosaic}
\begin{algorithmic}[1]
\Require a dimension $n$ and constraint set $\mathcal{C}$
\Ensure an $n$-mosaic satisfying $\mathcal{C}$
\For{$a = 1,\ldots,\texttt{max\_attempts}$}
    \State $M \gets$ an $n \times n$ blank matrix
    \For{$i = 0,\ldots,n-1$}
        \For{$j = 0,\ldots,n-1$}
            \State $A \gets \textsc{PotentialTiles}(M,i,j)$
            \State $M_{ij} \gets$ a uniformly random element of $A$
        \EndFor
    \EndFor
    \If{$M$ satisfies $\mathcal{C}$}
        \State \Return $M$
    \EndIf
\EndFor
\State raise an error
\end{algorithmic}
\end{algorithm}

For constrained generation, such as requiring a specified number of crossings or components, the function uses rejection sampling with a recursion depth limit of $5000$ attempts. This is effective for small mosaics; large mosaics with tight constraints may require multiple calls or more specialized search strategies.

The \texttt{is\_unknot} method uses the planar diagram interface to construct a knot object and compute knot Floer homology. Ozsváth and Szabó proved that knot Floer homology detects Seifert genus~\cite[Theorem~1.2]{ozsvathszabo2004genus}. Since the unknot is the unique knot of genus zero, and since the total rank of $\widehat{HFK}$ is one for the unknot, the method uses the criterion
\[
    \operatorname{rank} \widehat{HFK}(K) = 1
\]
to certify that a one-component mosaic represents the unknot. The corresponding routine is shown in Algorithm~\ref{alg:unknot-check}.

\begin{algorithm}[H]
\caption{Unknot detection using knot Floer homology}\label{alg:unknot-check}
\begin{algorithmic}[1]
\Require a suitably connected one-component mosaic $M$
\Ensure \texttt{True} if $\widehat{HFK}$ detects $M$ as the unknot
\If{$M$ has no crossings}
    \State \Return \texttt{True}
\EndIf
\State $K \gets \textsc{Link}(\textsc{PlanarDiagramCode}(M))$
\State $H \gets \textsc{KnotFloerHomology}(K)$
\State \Return $\operatorname{total\_rank}(H)=1$
\end{algorithmic}
\end{algorithm}

\subsection{Rational tangles}

The package also provides constructors for rational tangle mosaics. The $\infty$-tangle is represented by tile~$T_7$, the $0$-tangle by~$T_8$, and integer tangles by modified Jordan block matrices using crossing tiles $T_9$ for negative values and $T_{10}$ for positive values. The function \texttt{tangle\_join} arranges two rational tangle mosaics in a block matrix with connector tiles, providing a starting point for computational experiments with rational knots and links.

\section{Worked Examples}\label{sec:examples}

We now illustrate the main workflows of the package. The examples are intentionally small enough to read directly from the displayed matrices, except for the final random example, which demonstrates a larger generated mosaic.

\subsection{Instantiation and visualization}

We construct a $5 \times 5$ knot mosaic representing the five-crossing twist knot $5_2$ and verify its properties.

\begin{lstlisting}
(*@\sageprompt@*) M = Mosaic([[0,2,1,0,0],
(*@\sageprompt@*)             [2,9,10,1,0],
(*@\sageprompt@*)             [3,10,9,10,1],
(*@\sageprompt@*)             [0,3,7,8,4],
(*@\sageprompt@*)             [0,0,3,4,0]])
(*@\sageprompt@*) M.show()
(*@\qquad \includegraphics[width=0.25\linewidth]{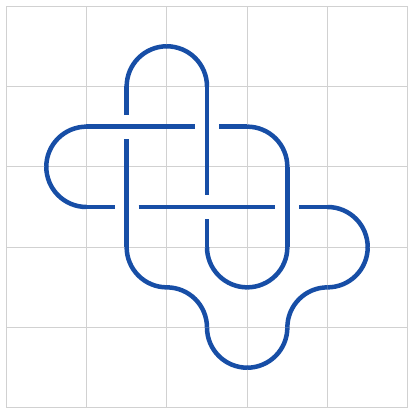}@*)
(*@\sageprompt@*) M.is_suitably_connected()
(*@\qquad \ @*) True
(*@\sageprompt@*) M.number_of_crossings()
(*@\qquad \ @*) 5
\end{lstlisting}

\subsection{Strand tracing and planar diagram codes}

The \texttt{strand\_of} method traces a strand from a given starting tile. The optional \texttt{verbose} parameter includes direction information in the output, returning a list of (coordinate, direction) pairs rather than coordinates alone:

\begin{lstlisting}
(*@\sageprompt@*) M.strand_of((1,1), direction='right', verbose=True)
(*@\qquad \ @*) Went right into tile (1, 2).
(*@\qquad \ @*) Went right into tile (1, 3).
(*@\qquad \ @*) Went down into tile (2, 3).
(*@\qquad \ @*) ... (24 steps total, returning to (1, 1))
\end{lstlisting}

\noindent\begin{minipage}{\linewidth}
\noindent The \texttt{pd\_code} method converts a mosaic to a planar diagram code, enabling computation of classical knot invariants:

\begin{lstlisting}
(*@\sageprompt@*) code = M.pd_code(); code
(*@\qquad \ @*) [[9, 4, 10, 5], [3, 10, 4, 1], [5, 8, 6, 9], [1, 6, 2, 7], [7, 2, 8, 3]]
(*@\sageprompt@*) K = Link(code)
(*@\sageprompt@*) K.jones_polynomial()
(*@\qquad \ @*) 1/t - 1/t^2 + 2/t^3 - 1/t^4 + 1/t^5 - 1/t^6
\end{lstlisting}
\end{minipage}

The same planar diagram bridge can be used to compare two different mosaic representatives. For example, if $M$ is the $5 \times 5$ mosaic above, the following mosaics of varying dimensions represent isotopic knots:

\begin{figure}[H]
\centering
\includegraphics[width=0.95\linewidth]{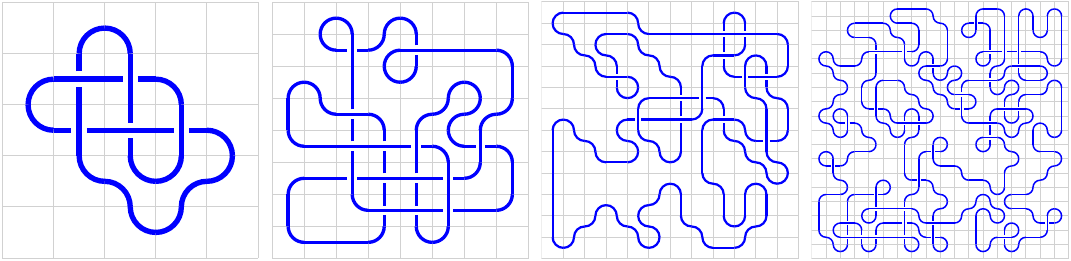}
\figcaptiongap
\caption{Four isotopic knot mosaic representations of the $5_2$ knot, as detected by the \texttt{is\_isotopic} function. By the Lomonaco--Kauffman conjecture, all these mosaics are related by a finite sequence of mosaic Reidemeister moves.}\label{fig:isotopic-mosaics}
\end{figure}

This can be confirmed by comparing the corresponding SageMath links for any pair of mosaics $M$ and $N$ under the \texttt{is\_isotopic} method:

\begin{lstlisting}
(*@\sageprompt@*) Link(M.pd_code()).is_isotopic(Link(N.pd_code()))
(*@\qquad \ @*) True
\end{lstlisting}




\subsection{Zooming and walking}

\begin{figure}[ht]
\centering
\includegraphics[width=0.6\linewidth]{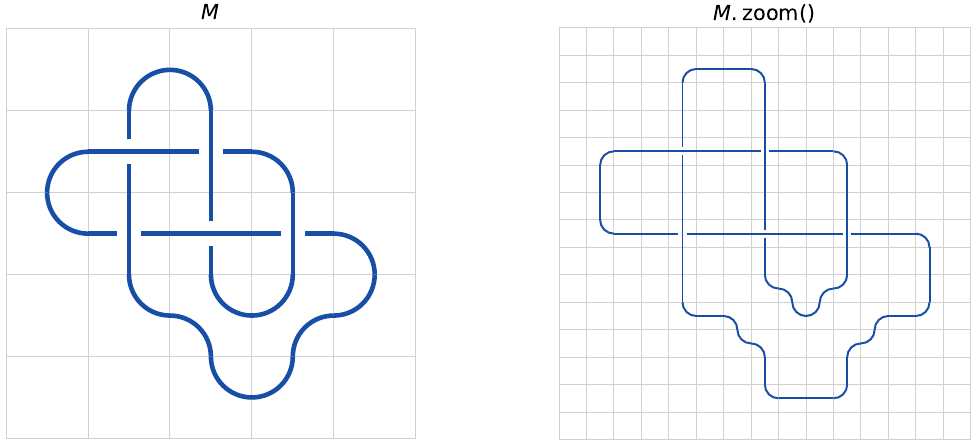}
\figcaptiongap
\caption{The zoom map applied to a mosaic, replacing each tile by an isotopy-equivalent $3 \times 3$ block.}\label{fig:zoom}
\end{figure}

The zoom operation replaces each tile by an equivalent $3 \times 3$ block. In the larger mosaic, each crossing is separated from its neighbors by non-crossing tiles. This is essential for the \texttt{walk} method, which traces paths between crossings. Figure~\ref{fig:zoom} illustrates the zoom map applied to a small mosaic.

\begin{lstlisting}
(*@\sageprompt@*) W = M.zoom()
(*@\sageprompt@*) crossings = W.find_crossings()
(*@\sageprompt@*) W.walk(crossings[0], 'right', path_list=True)
(*@\qquad \ @*) [(4, 4), (4, 5), (4, 6), (4, 7)]
\end{lstlisting}

\subsection{Rational tangles}

The package supports construction of rational tangle mosaics, following Conway's algebraic tangle notation~\cite{conway1970enumeration}. Positive integers produce mosaics using~$T_{10}$ crossings, negative integers use~$T_9$ crossings, and the special values~$0$ and~$\infty$ produce single-tile tangles.

\begin{figure}[H]
\centering
\begin{minipage}{0.48\linewidth}
\centering
\includegraphics[width=\linewidth]{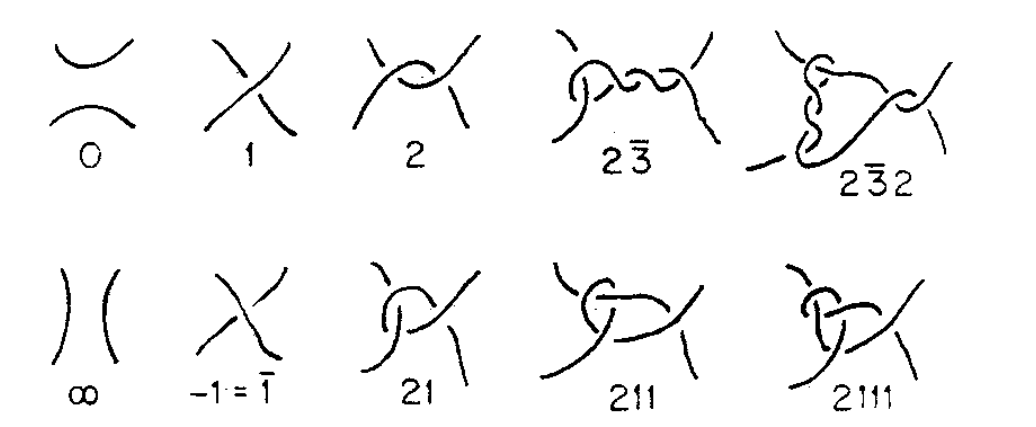}
\end{minipage}\hfill
\begin{minipage}{0.48\linewidth}
\centering
\includegraphics[width=\linewidth]{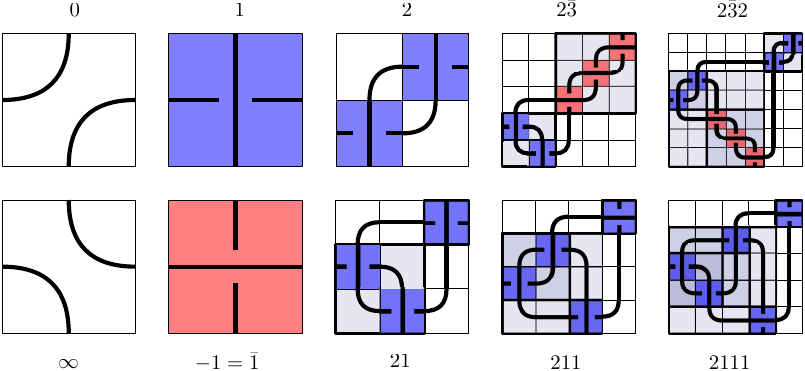}
\end{minipage}
\figcaptiongap
\caption{Left: Conway's original tangle examples from~\cite[Figure~1]{conway1970enumeration}. Right: corresponding mosaic rational tangles, including the $\infty$-tangle ($T_7$), $0$-tangle ($T_8$), and positive and negative integer tangles. The mosaic representation stores these diagrams as tile matrices.}\label{fig:tangles}
\end{figure}

The corresponding constructors can be called directly, returning mosaics that can be drawn or used in later computations:

\begin{lstlisting}
(*@\sageprompt@*) rational_tangle(-5).show()
(*@\qquad \includegraphics[width=0.18\linewidth]{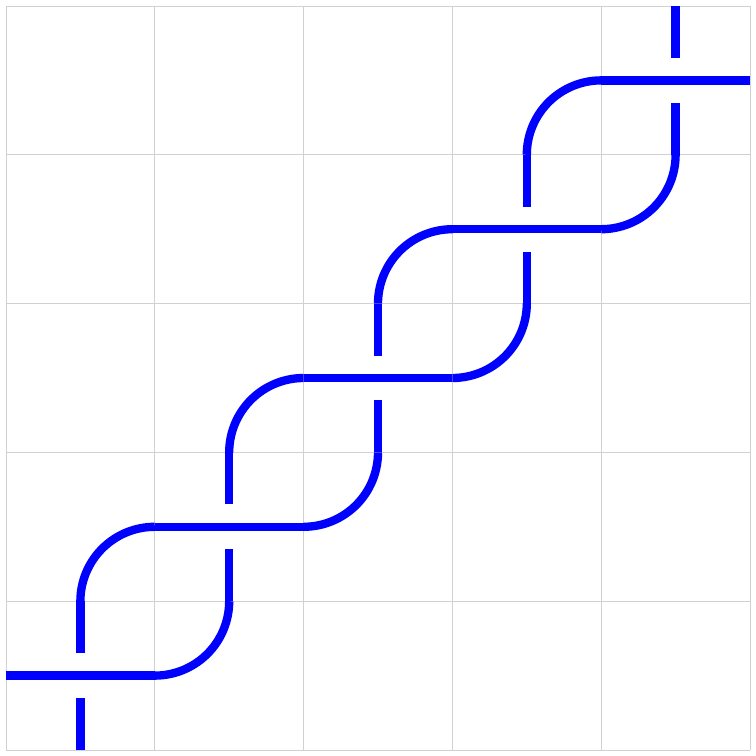}@*)
\end{lstlisting}

\noindent Figure~\ref{fig:tangles} compares Conway's original tangle notation with the corresponding mosaic-style tangle diagrams. The mosaic grid makes the recursive block structure explicit: each tangle is encoded as a finite tile matrix, so construction, visualization, and later computation all use the same combinatorial object.

Two tangles can be joined into a single mosaic:

\begin{lstlisting}
(*@\sageprompt@*) K = tangle_join([6,4])
(*@\sageprompt@*) K.show()
(*@\qquad \includegraphics[width=0.2\linewidth]{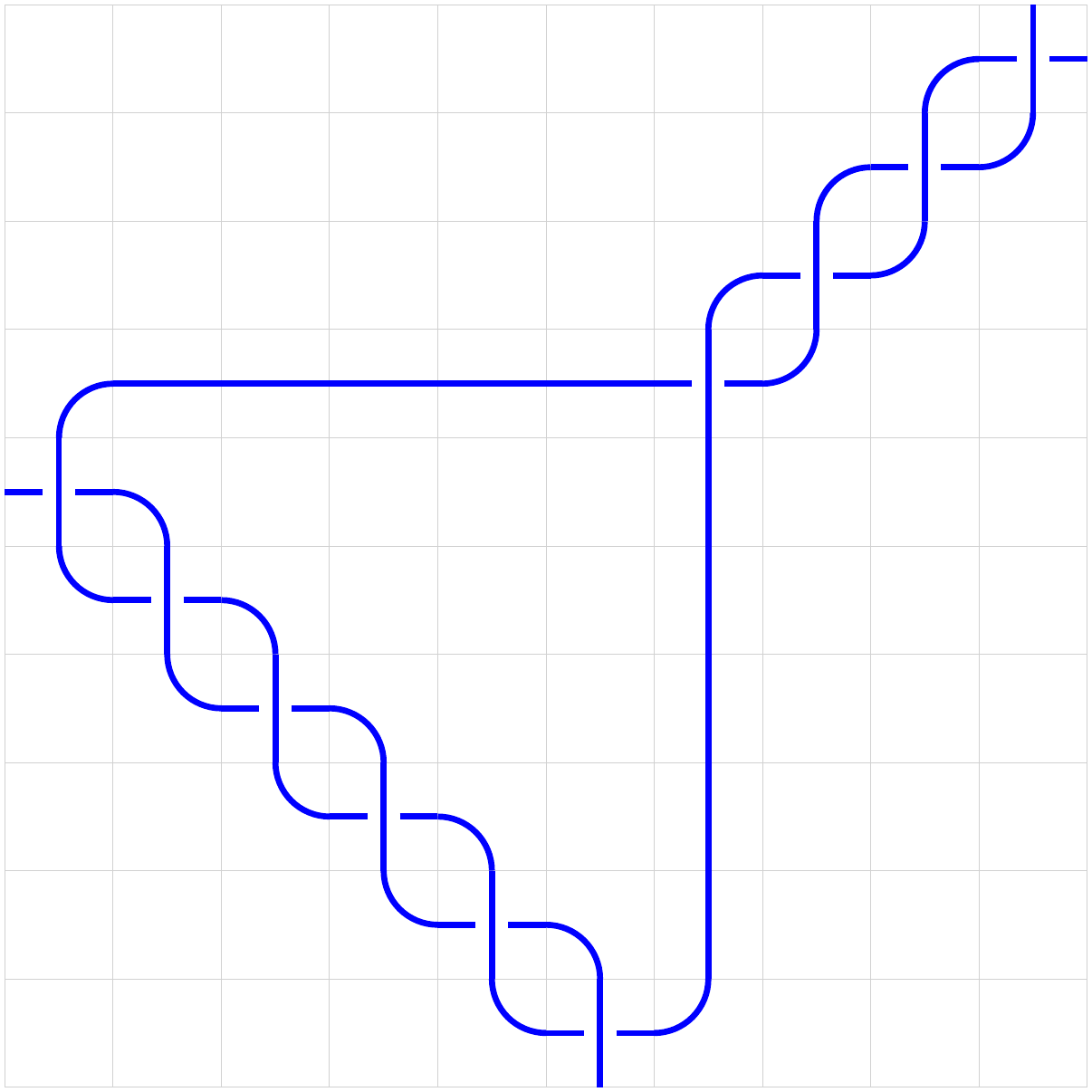}@*)
(*@\sageprompt@*) K.matrix()
(*@\qquad \ @*) [ 0  0  0  0  0  0  0  0  2 10]
(*@\qquad \ @*) [ 0  0  0  0  0  0  0  2 10  4]
(*@\qquad \ @*) [ 0  0  0  0  0  0  2 10  4  0]
(*@\qquad \ @*) [ 2  5  5  5  5  5 10  4  0  0]
(*@\qquad \ @*) [10  1  0  0  0  0  6  0  0  0]
(*@\qquad \ @*) [ 3 10  1  0  0  0  6  0  0  0]
(*@\qquad \ @*) [ 0  3 10  1  0  0  6  0  0  0]
(*@\qquad \ @*) [ 0  0  3 10  1  0  6  0  0  0]
(*@\qquad \ @*) [ 0  0  0  3 10  1  6  0  0  0]
(*@\qquad \ @*) [ 0  0  0  0  3 10  4  0  0  0]
\end{lstlisting}

\subsection{Multi-component mosaics}

The package handles both knots and links. For example, the following $4 \times 4$ mosaic represents a $2$-component link with $4$ crossings:

\begin{lstlisting}
(*@\sageprompt@*) H = Mosaic([[0,2,1,0],[2,9,10,1],[3,10,9,4],[0,3,4,0]])
(*@\sageprompt@*) H.show()
(*@\qquad \includegraphics[width=0.2\linewidth]{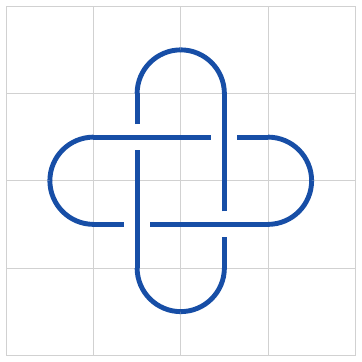}@*)
(*@\sageprompt@*) H.is_suitably_connected()
(*@\qquad \ @*) True
(*@\sageprompt@*) H.number_of_components()
(*@\qquad \ @*) 2
\end{lstlisting}

\noindent The planar diagram code interface records both components of this example, so it can be passed to link software in the same way as the one-component example above.

\subsection{Random mosaic generation}

Random suitably connected mosaics can be generated with optional constraints:

\begin{lstlisting}
(*@\sageprompt@*) R = random_mosaic(4)                                # any 4-mosaic
(*@\sageprompt@*) R = random_mosaic(4, number_of_crossings=2)         # exactly 2 crossings
(*@\sageprompt@*) R = random_mosaic(5, number_of_components=1)        # a knot (1 component)
(*@\sageprompt@*) R = random_mosaic(4, suitably_connected=False)      # any tiles
\end{lstlisting}

The \texttt{unknot} constraint filters the random search using \texttt{is\_unknot}. For example, the following code produces a $15 \times 15$ unknot mosaic.

\begin{lstlisting}
(*@\sageprompt@*) M = random_mosaic(15, unknot=True)
(*@\sageprompt@*) M.show()
(*@\qquad \includegraphics[width=0.45\linewidth]{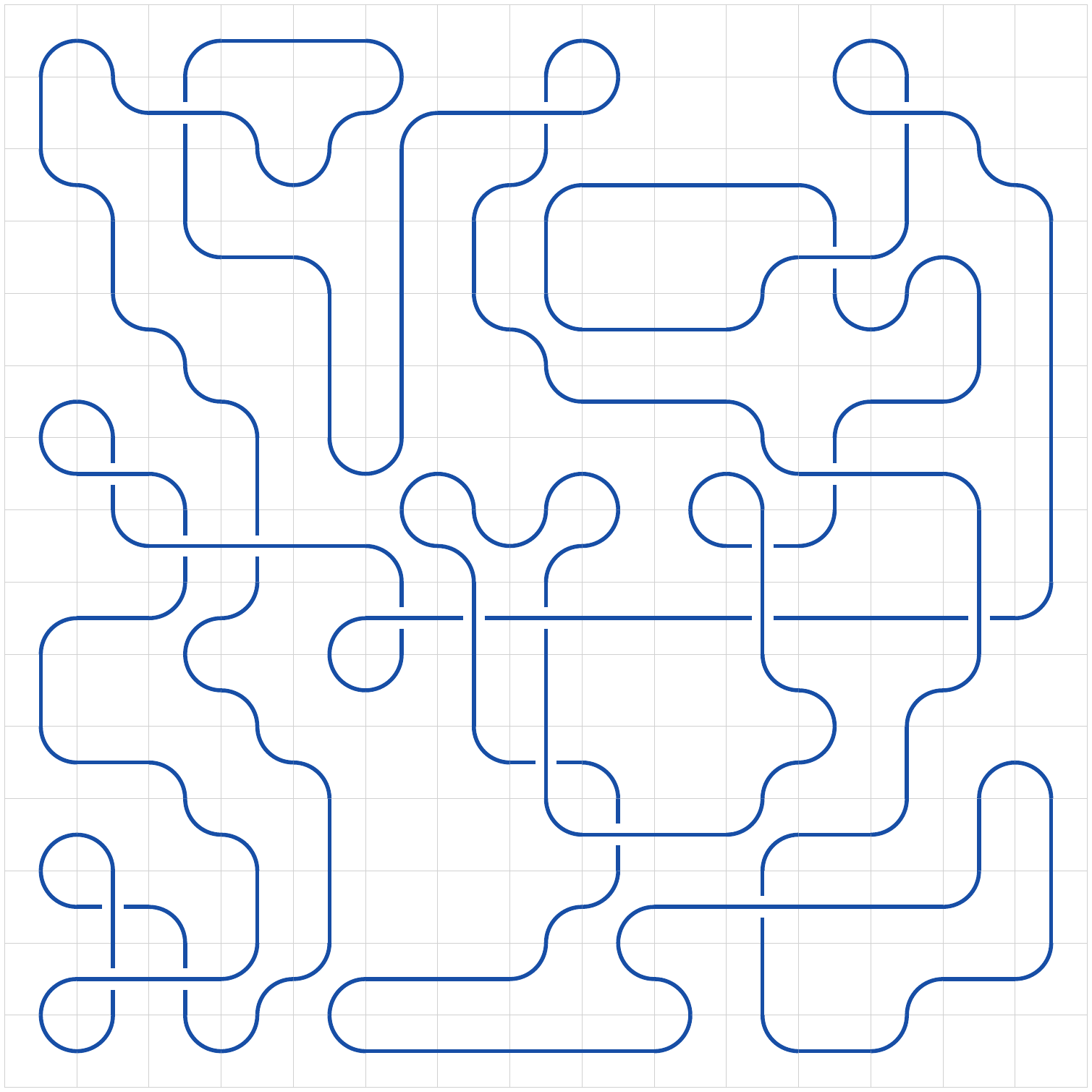}@*)
\end{lstlisting}

\section{Computational Applications and Future Work}\label{sec:applications}

The \texttt{KnotMosaics} package facilitates computational exploration of open questions in knot mosaic theory. We describe several directions where the current implementation is already useful and where further development would be valuable.

\subsection{Mosaic number and tile number}

The \emph{mosaic number} $m(K)$ of a knot~$K$ is the minimum~$n$ such that $K$ can be represented as a knot $n$-mosaic. Lee, Hong, Lee, and Oh~\cite{02d_Lee2014-3-5} proved that for any nontrivial knot or non-split link (other than the Hopf link), $m(K) \le c(K) + 1$, where $c(K)$ denotes the crossing number. Heap and Knowles~\cite{heapknowles2018} introduced the \emph{tile number} $t(K)$, the minimum number of non-blank tiles in any mosaic representation, and established bounds relating it to the mosaic number. Mosaic numbers and tile numbers have been tabulated for all prime knots with crossing number at most~$10$~\cite{geneseo2023}.



\subsection{Enumeration and growth}

Oh, Hong, Lee, and Lee~\cite{02d_Lee2014-4-5} established recursive formulas for the total number of knot $n$-mosaics using state matrix methods. Oh~\cite{oh2016growth} subsequently proved the existence of a \emph{knot mosaic growth constant}~$\delta$ satisfying $4 \le \delta \le (5 + \sqrt{13})/2$, where $\delta = \lim_{n \to \infty} D_n^{1/n^2}$ and $D_n$ denotes the total number of knot $n$-mosaics. The package's \texttt{potential\_tiles} method, which computes valid tile placements at each position, can serve as the basis for exhaustive enumeration algorithms to verify and extend these counts.

\subsection{Tangle operations and rational knots}

The \texttt{rational\_tangle} and \texttt{tangle\_join} functions allow systematic construction of rational knot mosaics from continued fraction data, providing a computational framework for investigating the relationship between tangle arithmetic and mosaic complexity. For instance, one can explore the mosaic number of rational knots as a function of their continued fraction expansion.

\subsection{Conclusion}

\texttt{KnotMosaics} makes knot mosaics available as explicit computational objects: matrices of tiles that can be checked, drawn, traversed, converted, and sampled. This direct connection between the combinatorial definition and executable code is useful both for producing examples and for testing conjectures about mosaic complexity. By exporting planar diagram codes, the package also connects mosaic-based computations to existing tools for classical knot and link invariants, allowing mosaic theory to be explored inside the broader \emph{SageMath} ecosystem.

\printbibliography

\end{document}